\documentclass[a4paper]{article}
\usepackage[utf8]{inputenc}
\usepackage{amssymb,amsmath,latexsym,mathtools}
\usepackage[english]{babel}

\usepackage{tikz-cd}

\title{Field Generated by Division Points of Certain Formal Group Laws}
\author{Soumyadip Sahu}
\date{}

\begin{document}
\maketitle
\begin{abstract}
In this article we study the Galois group of field generated by division points of special class of formal group laws and prove an equivalent condition for the group to be abelian. Further, we explore relations between the endomorphism ring of a formal group and the Galois group of field generated by division points. These results extend the results in \cite{the}, \cite{several}. 
\end{abstract}

\section{Introduction}
Let $p$ be an odd prime and let $K$ be a finite extension of $\mathbb{Q}_p$. Put $O_K$ to be the ring of integers of $K$, let $\mathfrak{p}_K$  
denote the unique maximal ideal of $O_K$ and $v_{\mathfrak{p}_K}(\cdot)$ be the valuation associated to it. Fix an algebraic closure $\overline{\mathbb{Q}}_p$ and $|.|_p$ be an fixed extension of the absolute value. Let $\overline{O}$ be the ring of integers of $\overline{\mathbb{Q}}_p$ and $\overline{\mathfrak{p}}$ be the unique maximal ideal of $\overline{O}$. Clearly $\overline{\mathfrak{p}} \cap K = \mathfrak{p}_K$.\\~\\
Let $A$ be an integrally closed, complete subring of $O_K$. Assume that $\mathfrak{p}_A$ is the maximal ideal, $K(A)$ is the field of fractions and $k(A)$ is the field of residues. Put $[k(A) : \mathbb{F}_p] = f_A$. Let $\mathfrak{F}$ be a (one dimensional, commutative) formal group-law defined over $O_K$ admitting an $A$ module structure. $\pi$ be a generator of $\mathfrak{p}_A$ and say \[[\pi](X) = \pi X + a_2X^2 + a_3X^3 + \cdots \in O_K[[X]] \tag{1.1}\] with at least one $a_i \in O_K - \mathfrak{p}_K$. Then $\min \,\{ i \,|\, |a_i|_p = 1 \} = p^h$ for some positive integer $h$ (see \cite{haz2}, 18.3.2). Now if $\pi_1$ is another generator of $\mathfrak{p}_A$ and \[[\pi_1](X) = \pi_1 X + b_2X^2 + \cdots \in O_K[[X]]\] then $b_{p^h} \in O_K - \mathfrak{p}_K$ and $\min \,\{ i \,|\, |b_i|_p = 1\} = p^h$. This integer $h$ is called the height of $\mathfrak{F}$ as formal $A$ module. If $a_i \in \mathfrak{p}_K$ for all $i \geq 2$ then we say, height of $\mathfrak{F}$ is infinity. We shall only consider formal $A$ modules of finite height.\\~\\
Let $A$ be as above and $\mathfrak{F}$ be a formal $A$ module over $O_K$. It defines a $A$ module structure on $\mathfrak{p}_K$ which naturally extends to a $A$-structure on $\overline{\mathfrak{p}}$. We shall denote the corresponding addition by $\oplus_{\mathfrak{F}}$ to distinguish it from usual addition.\\
Let $\pi$ be a generator of $\mathfrak{p}_A$. For each $n \geq 1$, use $\mathfrak{F}[\pi^n]$ to denote the $\pi^n$-torsion submodule of $\overline{\mathfrak{p}}$. For any sub-field $L$ of $\overline{\mathbb{Q}}_p$, let $L(\pi^n)$ be the subfield of $\overline{\mathbb{Q}}_p$ generated by $\mathfrak{F}[\pi^n]$ over $L$ and put $L(\pi^{\infty}) = \bigcup_{i \geq 1} L(\pi^i)$ We shall adopt the convention $\mathfrak{F}[\pi^0] = \{0\}$.\\~\\
Fix a generator $\pi$ of $\mathfrak{p}_K$ and put $K_{\pi} = \mathbb{Q}_p(\pi)$. Use $A_{\pi}$ to denote the ring of integers of $K_{\pi}$. For simplicity we shall write $\mathfrak{p}_{A_\mathfrak{\pi}}$ as $\mathfrak{p}_{\pi}$ and $f_{A_\pi}$ as $f_{\pi}$. Note that $\pi$ is a generator of $\mathfrak{p}_\pi$. Let $\mathfrak{F}$ be a formal $A_{\pi}$ module of height $h$ defined over $O_K$.\\~\\ 
The next remark summarizes main results proved in appendix-A of \cite{the}:\\~\\
\textbf{Remark 1.1 :} A)
i) With notation as above $f_{\pi}\,|\,h$.\\ 
Let $n \geq 1$. Put $q = p^h$ and $h_{r,\pi} = \frac{h}{f_{\pi}}$. Then the following statements are true :\\
ii) $\mathfrak{F}[\pi^n] \cong (A_{\pi}/\pi^nA_{\pi})^{h_{r,\pi}}$ as $A_{\pi}$ modules.\\
iii) If $z \in \mathfrak{F}[\pi^n] - \mathfrak{F}[\pi^{n-1}]$, then $K(z)/K$ is a totally ramified extension of degree $q^{n-1}(q - 1)$.\\
B) Results comparable to (A.i) and (A.ii)  already seem to exist in literature (for example, see \cite{lubin2}). Their complete proof is included in section-2 of appendix-A (\cite{the}).\\
C) In appendix-A of \cite{the} we made an assumption : \\
$\forall n \geq 1$, $K(z) = K(\pi^{n+1})$ for any $z \in \mathfrak{F}[\pi^{n+1}] - \mathfrak{F}[\pi^{n}]$.\\ We shall discuss more on this assumption in present article. Modulo this assumption the following theorem was proved :\\
Let $n, q, h_{r,\pi}$ be as before. In the following statements $h_{r,\pi}$ will be abbreviated as $h_r$.\\
Let $f$ be the degree of the extension of residue fields associated to the extension $K/\mathbb{Q}_p$. We shall assume that $h\,|\,f$.\\
Put $K^{h_{r}}_{\pi}$ to be the unique unramified extension of degree $h_r$ of $K_{\pi}$ in $\overline{\mathbb{Q}}_p$. Use $A_{\pi}^{h_r}$ to denote the ring of integers and $U_{i,K^{h_r}_{\pi}}$ ($i \geq 0$) to denote the $i$-th unit group of this field. Then :\\
i) $K(z) = K(\pi^n)$ for any $z \in \mathfrak{F}[\pi^n] - \mathfrak{F}[\pi^{n-1}]$.\\
ii) $K(\pi^n)/K$ is a Galois extension with \[\text{Gal}\,(K(\pi^n)|K) \cong U_{0,K_{\pi}^{h_r}}/U_{n,K_{\pi}^{h_r}}.\]
iii) Let $k$ and $i$ be integers such that $1 \leq k \leq n$ and $q^{k-1} \leq i \leq q^k - 1$. Then \[G_{i}(K(\pi^n)|K) = \text{Gal}(K(\pi^n)|K(\pi^{k})).\] 
D) The results of appendix-A of \cite{the} go through for $p=2$ also. \\~\\
\textbf{Remark 1.2 :} i) If $K/\mathbb{Q}_p$ is an unramified extension then one can take $\pi = p$. In this case $A_{\pi} = \mathbb{Z}_p$ and any formal group law over $O_K$ is $A_\pi$ module.\\
ii) Note that these results are similar to the results proved in classical Lubin-Tate theory (see \cite{neu}, chapter 3, section 6, 7 and 8). Though the hypothesis looks general it will turn out the assumption is quite strong and we shall explore equivalent versions of it in section-2 of this article.\\
iii) Consider the $\pi$-adic Tate module $T_{\pi}(\mathfrak{F}) = \varprojlim \mathfrak{F}[\pi^n]$. A priori it has $A_{\pi}$ module structure and the absolute Galois group $G_{K} = \text{Gal}(\overline{\mathbb{Q}}_p|K)$ acts on $T_{\pi}(\mathfrak{F})$ by $A_{\pi}$ module morphism. From remark-$1.1(A)$ it is easy to see that $T_{\pi}(\mathfrak{F})$ is a free $A_{\pi}$ module of rank $h_{r,\pi}$. But as a consequence of remark-1.1(C) (modulo the assumption) we can put a $A_{\pi}^{h_{r,\pi}}$ module structure on $T_{\pi}(\mathfrak{F})$ consistent with $A_{\pi}$ module structure. In this situation one can also see that $G_K$ acts on $T_{\pi}(\mathfrak{F})$ via $A_{\pi}^{h_{r,\pi}}$ module morphisms (Lemma- 3.18, lemma-3.19, appendix-A, \cite{the}).\\
iv) In appendix-A of \cite{the} results in section-2 and lemma 3.1-3.8 are independent of the assumption.\\
v) Let $A$ be a complete, integrally closed subring of $O_K$ containing $A_{\pi}$. If $\mathfrak{F}$ is a formal $A$ module of height $h$ defined over $O_K$ then we can prove analogues of results in remark-1.1 (A) and (C) with $h_{r,\pi}$ replaced by $\frac {h} {f_A}$. \\
vi) In this article we shall frequently refer to appendix-A and appendix-B without explicitly mentioning \cite{the}. \\
vii) The author feels that the results from appendix-A and appendix-B needs to studied more and currently the findings are being written up in parts (\cite{the}, \cite{several}). Once he finds the picture in satisfactory shape, he intends to rewrite them together in more comprehensive form. \\~\\
\textbf{Notations and conventions :}\\
 We often use notations introduced in this section without explicitly defining them again.\\ 
 Let $\pi$ be a generator of $\mathfrak{p}_K$ and let $\mathfrak{F}$ be a formal $A_{\pi}$ module over $O_K$. Such a formal module will be called a $\pi$-unramified group law over $O_K$.\\ 
 We shall use notations 
\[ \begin{split}
 q_h = p^h,\\
 h_{r,\pi} = \frac {h} {f_\pi}.
 \end{split}\] 
If $h, \pi$ are clear from context then we shall use the abbreviations $q, h_r$.\\
A $p$-adic field is a finite extension of $\mathbb{Q}_p$.\\
For $n \geq 1$, $\mu_{n}$ be the group of $n$-th roots of unity in $\overline{\mathbb{Q}}_p$.\\
If $L$ is a $p$-adic field then $O_L$ is the ring of integers and $\mathfrak{p}_L$ is the unique maximal ideal of $O_L$. Say residue degree of $L$ is $f$. Then $\mu_{p^f - 1} \subset O_L$ and $\mu_{p^f - 1} \cup \{0\}$ forms a cannonical set of representatives for the residue field. We shall always work with this set of representatives.\\~\\
\textbf{Acknowledgement :} I am thankful to Prof. S. David for pointing out an useful reference. I am thankful to Prof. C. Kaiser (Max Planck Institute for Mathematics , Bonn) for pointing out a mistake in an earlier version.  

\section{An equivalence}
In this section we study more about the assumption mentioned in  introduction. For this we need to introduce the language of $p$-divisible groups. We follow the treatment in \cite{tate}. Note that our set-up will be bit more general and we shall consider $\pi$-torsion points rather than $p$-torsion subgroups (if $K/\mathbb{Q}_p$ is unramified then one can choose $\pi = p$ ). The hypothesis and notation for present section is as in theorem-1.1 of introduction. \\~\\ 
Put $R = O_K[[X]]$. Define a $O_K$ algebra homomorphism $\phi : R \to R$ by defining $\phi(X) = [\pi](X)$. Note that $\phi$ makes $R$ into a $R$ module. Now we have the following lemma :\\~\\
\textbf{Lemma 2.1 :} The module structure defined above makes $R$ into a free $R$ module of finite rank.\\~\\
\textbf{Proof :} We shall show that $\{1, \cdots, X^{q-1}\}$ is a base .\\
First we shall show that the set is linearly independent. \\
Let $f_0(X), \cdots, f_{q-1}(X) \in R$ be such that \[\sum_{i = 0}^{q-1}\phi(f_i(X))X^i = 0.\tag{2.1}\]
We would like to show $f_i(X) = 0$ for each $0 \leq i \leq q-1$.\\
Assume that \[f_i(X) = \sum_{j \geq 0} a_{ij}X^j\] for each $0 \leq i \leq q-1.$\\
Put $m_i = \min \{v_{\mathfrak{p}_K}(a_{ij})\,|\,j \geq 1\}.$ for each $0 \leq i \leq q-1$ and $M = \min\,\{m_i \,|\, 0 \leq i \leq q-1\}$. If $M = \infty$ there is nothing to prove. Say, $M < \infty$. Define, $g_i(X) = \pi^{- M}f_i(X) \in R$ and $b_{ij} = \pi^{-M}a_{ij}$ for $0 \leq i \leq q-1$ and $j \geq 0$. From $(2.1)$ it follows that :\[\sum_{i = 0}^{q-1} \phi(g_i(X))X^i = 0.\]
Put $M_i = \min \{j \, |\, |a_{ij}|_p = 1\}$ and $m = \min \{ M_i \, |\, 0 \leq i \leq q-1 \}$. Clearly $m < \infty$. \\
Now modulo $\mathfrak{p}_K$ \[\sum_{i=0}^{q-1} \phi(g_i(X))X^{i} = u(b_{0m} + b_{1m}X + \cdots + b_{(q-1)m}X^{q-1})X^{mq} + \text{higher order terms}\] where $u$ is an unit and $b_{0m}, \cdots, b_{(q-1)m} \in O_K$.\\ 
By choice of $m$ at least one of $b_{im}$ is non-zero modulo $\mathfrak{p}_K$. Thus the sum must be non-zero modulo $\mathfrak{p}_K$ contrary to our assumption. Hence $M = \infty$ and each of $f_{i}(X)$ must be $0$ as desired.\\~\\
Now let $f(X) \in R$. We want to show that there are elements $a_0(X), \cdots,  a_{q-1}(X) \in R$ such that \[f(X) = \sum_{i = 0}^{q-1} \phi(a_i(X))X^i. \]
First we prove a claim :\\~\\
\textbf{Claim 2.2 :} There are elements $a_0, \cdots a_{q-1} \in O_K$ and $g(X) \in R$ such that \[f(X) = a_0 + \cdots + a_{q-1}X^{q-1} + [\pi](X)f_1(X).\]\\
\textbf{Proof :} First notice that there are elements $a_0(1), \cdots, a_{q-1}(1) \in O_K$ and $g(1)(X), f(1)(X) \in R$ such that \[ f(X) = a_0(1) + \cdots + a_{q-1}(1)X^{q-1} + [\pi](X)g(1)(X) + \pi (f(1)(X))\] since modulo $\mathfrak{p}_K$, $[\pi](X)$ is $X^q u_0(X)$ where $u_0(X) \in R$ is an unit. Now we can repeat the argument for $f(1)(X)$ and inductively construct elements $a_0, \cdots a_{q-1} \in O_K, f_1(X) \in R$ as desired. $\square$ \\~\\
We can repeat the argument for $f_1(X)$ and inductively construct elements $a_0(X), \cdots, a_{q-1}(X)$ satisfying \[f(X) = \phi(a_0(X)) + \cdots + \phi(a_{q-1}(X))X^{q-1}. \]
This completes the proof of lemma. $\square$\\~\\
Let $S = \text{Spf}\,(O_K[[X]])$. $\phi$ induces a $O_K$-map $S \to S$. Denote this map by $[\pi]$. Note that $\mathfrak{F}$ gives structure of a group-scheme on $S$ and $[\pi] : S \to S$ is a homomorhism with respect to this group structure. For $\nu \geq 0$ put $G_{\nu} = \text{Ker}([\pi^{\nu}])$. Note that $G_{\nu} = \text{Spec}(A_{\nu})$ where $A_{\nu} \cong R/[\pi^{\nu}](X)R$. So each $G_{\nu}$ is a connected finite group scheme over $O_K$ of order $p^{\nu h}$. We have an exact sequence: 
\[ 0 \xrightarrow{} G_{\nu} \xrightarrow{i_{\nu}} G_{\nu+1} \xrightarrow{[\pi^{\nu}]} G_{{\nu}+1} \] for each $\nu \geq 0$ where $i_{\nu}$ is the inclusion. Following Tate (\cite[section-2.2]{tate}) we can construct a $p$-divisible group $\mathfrak{F}_{\pi} = \varinjlim G_i$.\\~\\ 
Let $\text{Hom}_{O_K}(\mathfrak{F}, \mathfrak{F})$ be the set of all group homomorphisms of $\mathfrak{F}$ defined over $O_K$ and let $\text{Hom}_{O_K}^{A_\pi}(\mathfrak{F}, \mathfrak{F})$ denote the set of all $A_\pi$ module morphisms of $\mathfrak{F}$  defined over $O_K$. By a result of Hazewinkel any group endomorphism of $\mathfrak{F}$ is automatically an $A_\pi$ module morphism (see \cite[21.1.4]{haz2}). So these two sets are same and we shall use the notation $\text{Hom}_{O_K}(\mathfrak{F}, \mathfrak{F})$.\\~\\
Following the argument in \cite[section-2.2]{tate} we obtain a bijection $\text{Hom}_{O_K}(\mathfrak{F}, \mathfrak{F}) \to \text{Hom}_{O_K}(\mathfrak{F}_{\pi}, \mathfrak{F}_{\pi})$ induced by natural map. $T_{\pi}(\mathfrak{F})$ be the Tate module associated to $\mathfrak{F}_{\pi}$ (\cite[section-2.4]{tate}). Main result of the theory implies the canonical map $\text{Hom}_{O_K}(\mathfrak{F}_{\pi}, \mathfrak{F}_{\pi}) \to \text{Hom}_{G_K}(T_{\pi}(\mathfrak{F}), T_{\pi}(\mathfrak{F}))$ is bijective where $G_K = \text{Gal}(\overline{\mathbb{Q}}_p/K)$  and we are considering $T_{\pi}(\mathfrak{F})$ as $G_K$ module.  This proves that the natural map \[\text{Hom}_{O_K}(\mathfrak{F}, \mathfrak{F}) \to \text{Hom}_{G_K}(T_{\pi}(\mathfrak{F}), T_{\pi}(\mathfrak{F})) \tag{2.2}\] is bijective. \\~\\
With this set-up we have :\\~\\
\textbf{Theorem 2.3 :} Assume $h \,|\,f$ as in hypothesis of theorem-1.2. $A_{\pi}^{h_r}$ be as in remark-$1.3(iii)$. The following are equivalent :\\
i) $\text{Gal}(K(\pi^\infty)|K)$ is abelian.\\
ii) Assumption (from remark-1.1(C)) holds for each $n \geq 1$.\\
iii) $\mathfrak{F}$ has a formal $A_{\pi}^{h_r}$ module structure over $O_K$ which extends $A_{\pi}$ module structure.\\~\\
\textbf{Proof :} We shall show $(i) \iff (ii)$ and $(ii) \iff (iii)$.\\~\\
$(i) \implies (ii)$ We know that $K(\pi^{n+1})|K$ is Galois and given $z_1,  z_2 \in \mathfrak{F}[\pi^{n+1}] - \mathfrak{F}[\pi^n]$ there is a $\tau \in G_K$ such that $\tau(z_1) = z_2$ (see lemma-3.1 and lemma-3.8 in appendix-A). It follows that if $z \in \mathfrak{F}[\pi^{n+1}] - \mathfrak{F}[\pi^n]$ then $K(\pi^{n+1})$ is Galois closure of $K(z)$ over $K$. Each sub-extension of an abelian  extension is Galois. Hence the implication.\\ 
$(ii) \implies (i)$ Remark-$1.1(C)$ (in particular lemma-3.18, lemma-3.19 from appendix-A) implies \[\text{Gal}(K(\pi^\infty)|K) \cong U_{0,K_{\pi}^{h_r}} .\] So $\text{Gal}(K(\pi^\infty)|K)$ is abelian. Hence the implication.\\ 
$(ii) \implies (iii)\,$ Consider $\tau_{\infty} \in \text{Hom}_{G_K}(T_{\pi}(\mathfrak{F}),  T_{\pi}(\mathfrak{F}))$ (see remark-3.14 of appendix-A). Let $\tau_{\infty}(X) \in O_K[[X]]$ be the corresponding element in $\text{Hom}_{O_K}(\mathfrak{F}, \mathfrak{F})$. Since $\tau_{\infty}^{q-1} = \text{Id}$, $\tau_{\infty}^{q-1}(X) = X$ where $\tau^{q-1}$ denotes $(q-1)$-fold composition. From construction of the element $\tau_{\infty}$ (see lemma-3.13 in appendix-A) and the `analytic interpretation lemma' (lemma-3.5 in appendix-A) we conclude that there is a primitive $(q-1)$-th root of unity $\zeta$ such that \[\tau_{\infty}(X) = \zeta X + \text{higher degree terms}. \tag{2.3}\]
Note that $A_{\pi}^{h_r} = A_{\pi}[\zeta] = A_{\pi} \oplus A_{\pi}\zeta \oplus \cdots \oplus A_{\pi}\zeta^{h_r - 1}$ as $A_{\pi}$ module.\\
Thus one can define a ring homomorphism \[ [\cdot] : A_{\pi}^{h_r} \to \text{Hom}_{O_K}(\mathfrak{F}, \mathfrak{F}) \] by $\zeta \to \tau_{\infty}(X)$ and extending $A_{\pi}$ linearly. This is well defined due to the direct sum decomposition mentioned above. It is easy to check that for each $\alpha \in A_{\pi}^{h_r}$ \[[\alpha](X) = \alpha X + \text{higher degree terms} .\]
So we have the desired $A_{\pi}^{h_r}$ module structure.\\
$(iii) \implies (ii)$ First observe that $[k(A_{\pi}^{h_r}) : \mathbb{F}_p] = h$. From remark-$1.3(v)$ $\mathfrak{F}[\pi^{n+1}]$ is a free $A_{\pi}^{h_r}/\pi^{n+1}A_{\pi}^{h_r}$ module of rank $1$ for all $n \geq 1$. Let $z \in \mathfrak{F}[\pi^{n+1}] - \mathfrak{F}[\pi^{n}]$.  Then $z$ generates $\mathfrak{F}[\pi^{n+1}]$ as $A_{\pi}^{h_r}$ module. So $\mathfrak{F}[\pi^{n+1}] \subseteq K(z)$. Hence assumption is verified.\\
This completes the proof. $\square$ \\~\\
\textbf{Remark 2.4 :} A) $(iii) \implies (i)$ is a well-known result in literature (\cite{Serre}).\\
B) Let $L$ be an unramified extension of $K$ (possibly infinite). Use $\widehat{L}$ denote completion of $L$, $O_{\widehat{L}}$ for the ring of integers of $\widehat{L}$ and $\mathfrak{p}_{\widehat{L}}$ be the maximal ideal of $O_{\widehat{L}}$. If $\pi$ is a generator of $\mathfrak{p}_K$, then $\pi$ also generates $\mathfrak{p}_{\widehat{L}}$. Let $\mathfrak{F}$ be a $A_{\pi}$ module of height $h$ defined over $O_{\widehat{L}}$. Then one can prove analogues of results in appendix-A (see remark-1.1) for $\mathfrak{F}$ (note that in this case we shall replace the fields $K(\pi^n)$ by $\widehat{L}(\pi^n)$ and $K$ by $\widehat{L}$.) Proof is similar to proof in appendix-A. One can also prove an analogue of theorem-2.3.\\~\\

\section{Further results} 
Theorem-2.3 establishes an equivalence between properties of two different objects namely $\text{Gal}(K(\pi^{\infty})|K)$ and $\text{Hom}_{O_K}(\mathfrak{F}, \mathfrak{F})$. In this section we study them separately.\\
For convenience we shall denote $\text{Hom}_{O_K}(\mathfrak{F}, \mathfrak{F})$ by $\text{End}_{O_K}(\mathfrak{F})$. This object has been extensively studied by Lubin (see \cite{lubin1}, \cite{lubin2}). In first part of this section we recall some results from these two articles and rephrase some of them in context of $A_{\pi}$ modules. Remember that any group endomorphism of $\mathfrak{F}$ is an $A_{\pi}$ module morphism.
\subsection{The endomorphism ring} 
Let $\mathfrak{F}$ be a $A_{\pi}$ module of finite height defined over $O_K$ for some generator $\pi$ of $\mathfrak{p}_K$.\\
Consider the map \[c : \text{End}_{O_K}(\mathfrak{F}) \to O_K\]
defined by $c(f(X)) = a_1$ where $f(X) = \sum_{i \geq 1} a_iX^i \in O_K[[X]]$.\\
By lemma-2.1.1 in \cite{lubin1} $c$ is an isomorphism onto a closed subring of $O_K$. From the hypothesis of $A_\pi$ module structure on $\mathfrak{F}$, it follows that $A_\pi \subseteq c(\text{End}_{O_K}(\mathfrak{F}))$.\\
Note that $\text{End}_{O_K}(\mathfrak{F})$ is noethrian, local ring. It contains a subring isomorphic to $A_{\pi}$ which by abuse of notation will be identified with $A_{\pi}$. \\~\\
Let $T_{\pi}(\mathfrak{F})$ be the $\pi$-adic Tate module. It is a free $A_{\pi}$ module of rank $h_r$. Put, 
\[
\begin{split}
V_{\pi}(\mathfrak{F}) = T_{\pi}(\mathfrak{F}) \otimes_{A_\pi} K_{\pi}\\
\Lambda_{\pi}(\mathfrak{F}) = \bigcup_{i \geq 1} \mathfrak{F}[\pi^i].
\end{split}
\] 
We have an isomorphism of $A_{\pi}$ modules \[ \lambda : V_{\pi}(\mathfrak{F})/T_{\pi}(\mathfrak{F}) \to \Lambda_{\pi}(\mathfrak{F}).\]
Let $L$ be a free $A_{\pi}$ module of rank $h_r$ contained in $V_{\pi}(\mathfrak{F})$ (will be called a `full lattice') which is stable under the action of absolute Galois group $G_K$. Following the argument in \cite{lubin2}, with $L$ one can associate an $A_{\pi}$ module $\mathfrak{G}$ defined over $O_K$ which is isogenus to $\mathfrak{F}$ and this association has certain functorial properties (see theorem-2.2, \cite{lubin2}). Further if $L =  T_{\pi}(\mathfrak{F})$, then one can take $\mathfrak{G} = \mathfrak{F}$ and and the isogeny to be identity.\\
The following lemma is a modified version of theorem-3.1 from \cite{lubin2} :\\~\\
\textbf{Lemma 3.1.1 :} Let $\mathfrak{F}$ be a $A_{\pi}$ module of finite height over $O_K$ and assume that $L$ is $G_K$ stable full lattice in $V_{\pi}(\mathfrak{F})$ giving rise to a $A_{\pi}$ module $\mathfrak{G}$ over $O_K$ as described above. Note that since $\mathfrak{F}$, $\mathfrak{G}$ are isogenous, $c(\text{End}_{O_K}(\mathfrak{F}))$ and $c(\text{End}_{O_K}(\mathfrak{G}))$ have same fraction field in $K$ (to be denoted $F$). Let $\zeta \in F$. Now $\zeta \in \text{End}_{O_K}(\mathfrak{G})$ if only if $\zeta L \subseteq L$.\\~\\
\textbf{Proof :} Proof is similar to proof of theorem-3.1 from \cite{lubin2}. $\square$ \\~\\
We shall use lemma-3.1.1 to prove the following result which is a generalised version of theorem-3.3 from \cite{lubin2} and important for later development :\\~\\
 \textbf{Theorem 3.1.2 :} Let $\mathfrak{F}$ be an $A_{\pi}$ module of finite height over $O_K$ for some generator $\pi$ of $\mathfrak{p}_K$. Then $\text{End}_{O_K}(\mathfrak{F})$ is integrally closed in its fraction field.\\~\\
 \textbf{Proof :} From lemma-3.8 in appendix-A we know that the absolute Galois group $G_K$ acts transitively on $\mathfrak{F}[\pi] - \{0\}$. Hence only proper $G_K$ submodule of $\mathfrak{F}[\pi]$ is $\{0\}$. We have an exact sequence of $G_K$ modules :
 \[0 \to [\pi]T_{\pi}(\mathfrak{F}) \to T_{\pi}(\mathfrak{F}) \to \mathfrak{F}[\pi] \to 0.\]
From the observation above it follows that if $L$ is a $G_{K}$ stable full lattice containing $[\pi]T_{\pi}(\mathfrak{F})$ and contained in $T_{\pi}(\mathfrak{F})$ then either $L = [\pi]T_{\pi}(\mathfrak{F})$ or $L = T_{\pi}(\mathfrak{F})$.\\
Let $\mathfrak{m}$ be the maximal ideal of $\text{End}_{O_K}(\mathfrak{F})$. Clearly $[\pi] \in \mathfrak{m}$. Put $L = \mathfrak{m}T_{\pi}(\mathfrak{F})$. Clearly $L$ is a $G_K$ stable $A_{\pi}$ submodule of $T_{\pi}(\mathfrak{F})$ containing $[\pi]T_{\pi}(\mathfrak{F})$. From structure theorem of finitely generated modules over PID it follows that $L$ is free over $A_{\pi}$ of rank $h_r$. By the observation above $L = [\pi]T_{\pi}(\mathfrak{F})$ or $L = T_{\pi}(\mathfrak{F})$. Since $\text{End}_{O_K}(\mathfrak{F})$ is noetherian and $T_{\pi}(\mathfrak{F})$ is a finitely generated non-zero module it is clear that $L \neq T_{\pi}(\mathfrak{F})$. So $L = [\pi]T_{\pi}(\mathfrak{F})$.\\
Let $\phi \in \mathfrak{m}$. We shall show that there is a $\phi_1 \in \text{End}_{O_K}(\mathfrak{F})$ such that $\phi = [\pi]\phi_1$. Consider the element $\zeta = c(\phi)\pi^{-1}$ in fraction field of $c(\text{End}_{O_K}(\mathfrak{F}))$. From the equality above it follows that $\zeta T_{\pi}(\mathfrak{F}) \subseteq T_{\pi}(\mathfrak{F})$. By lemma-3.1.1 there is a $\phi_1 \in \text{End}_{O_K}(\mathfrak{F})$ such that $c(\phi_1) = \zeta$. This $\phi_1$ satisfies $\phi = \pi \phi_1$.\\
The argument above proves that $\mathfrak{m}$ is a principal ideal generated by $\pi$. Since $\text{End}_{O_K}(\mathfrak{F})$ is a local, noetherian domain we conclude that it is a DVR and hence the theorem. $\square$\\~\\
\textbf{Remark 3.1.3 :} $\text{End}_{O_K}(\mathfrak{F})$ is isomorphic to a closed subring of $O_K$ containing $A_{\pi}$. Let $F$ be the fraction field of $c(\text{End}_{O_K}(\mathfrak{F}))$. Then we have a tower of fields $K_{\pi} \subseteq F \subseteq K$. Since $K/K_{\pi}$ is unramified, $F/K_{\pi}$ is also unramified. The theorem above shows that $c(\text{End}_{O_K}(\mathfrak{F})) = O_F$, the ring of integers of $F$. Thus $\mathfrak{F}$ has a $O_F$ module structure over $O_K$ extending $A_{\pi}$ module structure. Let $f_{O_F}$ be the degree of residue extension of $O_F$. Since $F/K_{\pi}$ is unramified $[F : K_{\pi}] = \frac {f_{O_F}} {f_\pi}$ and $\mathfrak{F}$ has same height $h$ with respect to $O_F$. From theorem-1.1 we have $f_{O_F} \,|\,h$. Hence $F \subseteq K_{\pi}^{h_r}$. Thus $K_{\pi} \subseteq F \subseteq K_{\pi}^{h_r}$. If $F = K_{\pi}^{h_r}$, we say the endomorphism ring of $\mathfrak{F}$ has full height.

\subsection{The Galois Group}
Let $\mathfrak{F}$ be a $\pi$-unramified group law of height $h$ defined over $O_K$. In this sub-section we study the group $\text{Gal}(K(\pi^{\infty})|K)$. For simplicity we shall use the notation $G_K^{\infty}$.\\
Note that $G_K$ acts $T_{\pi}(\mathfrak{F})$ by $A_{\pi}$ linear maps. $T_{\pi}(\mathfrak{F})$ is a free module $A_{\pi}$ of rank $h_r$. Let $\{z_1,\cdots, z_{h_r}\}$ be a fixed $A_{\pi}$ base for $T_{\pi}(\mathfrak{F})$. This choice identifies $\text{End}_{A_\pi}(T_{\pi}(\mathfrak{F})$ with $M_{h_r}(A_\pi)$ and invertible endomorphisms with $\text{Gl}_{h_r}(A_\pi)$.  Thus we  have a group homomorphism \[\rho(\mathfrak{F}) : G_{K} \to \text{Gl}_{h_r}(A_{\pi}).\]
Put $H = \rho(\mathfrak{F})(G_K)$.\\
Note that kernel of $\rho(\mathfrak{F})$ is $\text{Gal}(\overline{\mathbb{Q}}_p|K(\pi^{\infty}))$ and $H \cong G_K^{\infty}$ as an abstract group.\\ 
Further it is easy to see that $\rho(\mathfrak{F})$ is continuous with respect to usual profinite topology. Since $G_K$ is compact with respect to pro-finite topology, $H$ is closed subgroup of $\text{Gl}_{h_r}(A_{\pi})$ and the isomorphism of abstract groups mentioned above is actually an isomorphism of topological groups.\\
Let $R$ be the $A_{\pi}$ submodule of $M_{h_r}(A_{\pi})$ generated by $H$. Note that, we have a canonical surjective ring homomorphism \[\rho(\mathfrak{F}) : A_{\pi}[G_K] \to R \] induced by $\rho(\mathfrak{F})$ satisfying $a.z = \rho(\mathfrak{F})(a).z$ for all $a \in A_{\pi}[G_K]$ and $z \in T_{\pi}(\mathfrak{F})$.\\~\\
Let \[\phi : \text{End}_{O_K}(A_{\pi}) \to M_{h_r}(A_{\pi})\] be the canonical ring homomorphism. Put $\text{Im}(\phi) = E$.\\
Let $R'$ denote the commutant (ie. centralizer) of $R$ in $M_{h_r}(A_{\pi})$. Clearly $E \subseteq R'$. The bijection in $(2.2)$ proves that $\phi$ is injective and $E = R'$.\\~\\
$F$ be the fraction field of $c(\text{End}_{O_K}(A_\pi))$. Note that we have a tower of fields $K_{\pi} \subseteq F \subseteq K_{\pi}^{h_r}$ and $[F : K_{\pi}^{h_r}] = \frac {f_{F}} {f_{\pi}}$ where we are using the abbreviation $f_{F}$ for $f_{O_F}$. So, $O_{F} = A_{\pi}[\zeta]$ where $\zeta$ is a primitive $(p^{f_F} - 1)$-th root of unity and \[O_F = A_{\pi} \oplus \cdots \oplus A_{\pi}\zeta^{\frac {f_F} {f_\pi} - 1}\] as $A_{\pi}$ modules.\\ 
Note that $T_{\pi}(A_{\pi})$ is a free $O_F$ module of rank $h/f_F$. Let $\{z_1,\cdots, z_{h/f_F}\}$ be a $O_F$ base. Then $\{[\zeta]^iz_j \,|\, 0 \leq i \leq f_F/f_\pi -1, 1 \leq j \leq h/f_F \}$ is a $A_{\pi}$ base for $T_{\pi}(\mathfrak{F})$. We shall use such base now on, order it using reverse dictionary order (ie $(i,j) \leq (i_1, j_1)$ iff $(j,i) \leq (j_1, i_1)$ in dictionary order)  and assume that $\rho(\mathfrak{F})$ and $\phi$ are defined with respect to this ordered base. \\
In the following calculation we shall partition a $h_r \times h_r$ matrix into square blocks of $f_F/f_{\pi}$ matrices. For convenience put $m = f_F/f_{\pi}$ and $n = h/f_F$.\\
Now $\phi([\zeta]) = (A_{ij})_{n \times n}$ where $A_{ij} = 0_{m \times m}$ whenever $i \neq j$ and $A_{ii} = A$ for all $1 \leq i \leq n$ where
\[A =
\begin{bmatrix}
0 & 0 & . & 0 & 1\\
1 & 0 & . & 0 & 0\\
0 & 1 & . & 0 & 0\\
. & . & . & . & .\\ 
0 & 0 & . & 1&  0
\end{bmatrix}_{m \times m}\]
Let $Y = (y_{ij})_{m \times m}$ be a matrix whose entries are indeterminates. Note that the relation $AY = YA$ holds if and only if relations $\{R(i)\,|\, 0 \leq i \leq m-1\}$ hold, where $R(i)$ is defined as \[ R(i) : y_{1,1+i} = y_{2,2+i} = \cdots = y_{m,m+i} \] with the convention that if some index exceeds $m$ we shall read it modulo $m$.\\
Now $X = (X_{ij})_{n \times n}$ be a matrix in $h_r^2$ indeterminates where $X_{ij}$ is a $m \times m$ matrix with $X_{ij} = (x_{pq,ij})_{m \times m}$. Note that $X\phi([\zeta]) = \phi([\zeta])X$ if and only if $AX_{ij} = X_{ij}A$ for all $1 \leq i, j \leq n$. \\
Let $V$ be the sub-variety of $M_{h_r}(A_{\pi})$ defined by \[ V = \{ X \in M_{h_r}(A_\pi) \,|\, X\phi([\zeta]) = \phi([\zeta])X \}.\]
Clearly $H \subseteq R \subseteq V$.\\
From the relations above it is easy to see $V \otimes_{A_\pi} K_{\pi}$ has dimension $ \leq n^{2}m = \frac {h_r h} {f_{F}} $. Hence $H$ (thought as a subset of $M_{h_r}(K_\pi)$) is contained in a variety of dimension $\leq \frac {h_r h} {f_F}$.\\~\\
\textbf{Remark 3.2.1 :} i) Consider the case $f_F = h$. From theorem-2.3 and remark-1.1(C) it follows that $R \cong A_{\pi}^{h_r}$ and $H \cong U_{0,K_{\pi}^{h_r}}$ the unit group of this ring. Since $R = A_{\pi}[\zeta]$, it is a free $A_{\pi}$ module of rank $h_r$. It follows that $H$ is an algebraic subgroup of $\text{Gl}_{h_r}(A_{\pi})$ and $H \otimes_{A_\pi} K_{\pi}$ has dimension exactly $h_r$. This leads to the following questions about general situation :\\
a) Is $H$ is always an algebraic subgroup ? \\
b) Can $H$ (thought as subset of $G_{h_r}(K_{\pi})$) can be put inside a sub-variety of $\text{Gl}_{h_r}(K_\pi)$ of dimension $< \frac {h_r h} {f_F}$ ? \\
Further if $f_F = f_\pi$, is $H$ an open subgroup of $G_{h_r}(A_{\pi})$ (in usual profinite topology) ?\\
ii) Consider the reduction modulo $\pi$ map $\text{Gl}_{h_r}(A_\pi) \to \text{Gl}_{h_r}(\mathbb{F}_{p^{f_\pi}})$ and assume that $\widetilde{H}$ is the image of $H$. The kernel of $H \to \widetilde{H}$ is $\text{Gal}(K(\pi^{\infty})|K(\pi))$ and hence $\widetilde{H} \cong \text{Gal}(K(\pi)|K) \cong \mathbb{Z}/(q-1)\mathbb{Z}$. \\
iii) If $G_{K}^{\infty}$ is abelian then $H$ contains all the scalar matrices. Conversely, one can ask what can be said about $G_{K}^{\infty}$ or in particular about $f_F$ if $H$ contains all scalar matrices.\\
iv) Note that one can construct $A_{\pi}$ modules with distinct $f_{F}$ (subject to the conditions $f_{\pi} \,|\, f_{F} \,|\, h$). For example a Lubin-Tate module has $f_F = h$. By theorem-5.1.2 in \cite{lubin1} one can construct for any such $f_{F}$ if $K/\mathbb{Q}_p$ is unramified (with $\pi= p$).   
 
\section{Applications to theory of local fields}
This is the concluding section and its goal is to apply the theory developed in \cite{the}, \cite{several} and in this article  to construct extensions of local fields. This aspect requires more careful study and the author hopes to do so in some future article. \\~\\
We begin by some simple observations :\\~\\
\textbf{Remark 4.1 :} i) Lubin-Tate theory is used to generate totally ramified abelian extensions. If $\mathfrak{F}$ is an unramified group law of height $h$ over $O_K$ such that the endomorphism ring has full height, then $K(\pi^{\infty})/K$ is also a totally ramified abelian extension whose Galois group is isomorphic to the Galois group of Lubin-Tate extension of height $h$ (remark-1.1, theorem-2.3). Now let $\mathfrak{F}_1, \mathfrak{F}_2$ be two unramified group laws (wrt $\pi_1$ and $\pi_2$ respectively) of height $h$ over $O_K$ having endomorphism rings of full height. Clearly $K(\Lambda_{\pi_1}(\mathfrak{F}_1), \Lambda_{\pi_2}(\mathfrak{F}_2))$ is an abelian extension of $K$. It is interesting to ask what is the Galois group of this extension. Put $L = K^{\text{ur}}$. Let $\widehat{L}$ be the completion of $L$. Clearly $\text{Gal}(\widehat{L}(\Lambda_{\pi_1}(\mathfrak{F}_1))|\widehat{L}) \cong \text{Gal}(K(\Lambda_{\pi_1}(\mathfrak{F}_1))|K)$. One can consider the extension $\widehat{L}(\Lambda_{\pi_1}(\mathfrak{F}_1), \Lambda_{\pi_2}(\mathfrak{F}_2))/ \widehat{L}$ and ask for its Galois group. Note that if $\mathfrak{F}_1$ and $\mathfrak{F}_2$ are Lubin-Tate extensions then $\widehat{L}(\Lambda_{\pi_1}(\mathfrak{F}_1), \Lambda_{\pi_2}(\mathfrak{F}_2)) = \widehat{L}(\Lambda_{\pi_1}(\mathfrak{F}_1))$. Further, for any two unramified group laws of same height $\mathfrak{F}_1, \mathfrak{F}_2$ we have $\widehat{L}(\mathfrak{F}_1[\pi_1], \mathfrak{F}_2[\pi_2]) = \widehat{L}(\mathfrak{F}_1[\pi_1])$ (see \cite{several}, theorem-1.1). \\
ii) Same question can be asked if we have two unramified group laws of distinct height having endomorphism rings of full height.\\
iii) Now let $\mathfrak{F}$ be an unramified group law of height $h$ defined over $O_K$ such that the endomorphism ring does not have full height. Then $K(\pi^\infty)|K$ is a nonabelian extension. This gives a method of constructing non-abelian extensions explicitly. One can ask for a family of unramified formal group laws so that any non-abelian extension of $K$ is contained in one such of extension or inside compositum of several ones of them.\\~\\
We end the section with a result on roots of unity.

\subsection{Roots of unity}
Let $\mathfrak{F}$ be a unramified group law defined over $O_K$ corresponding to a generator $\pi$. It is natural to ask the question if $\zeta_{p^n} \in K(\pi^n)$ for some (all) $n \geq 1$, where $\zeta_{p^n}$ is a primitive $p^n$-th root of unity. We may relax the condition a bit and ask an easier question if $\zeta_{p^n} \in K^{\text{ur}}(\pi^n)$ where $K^{\text{ur}}$ is the maximal unramified extension of $K$ inside $\overline{\mathbb{Q}}_p$.\\
This question is easier to answer for some special group laws. Let $K$ be an unramified extension of $\mathbb{Q}_p$. If $\mathfrak{F}$ is the formal group law associated to a super-singular elliptic curve $E$ defined over $K$, then by Weil pairing $\zeta_{p^n} \in K(p^n)$ for all $n \geq 1$ (Note that in this case $E[p^n] = \mathfrak{F}[p^n]$). Further if $\mathfrak{F}$ is the formal group law associated to an elliptic curve with good, ordinary reduction then $E[p^n] \subseteq K^{\text{ur}}(p^n)$ for all $n \geq 1$ and by Weil pairing we have $\zeta_{p^n} \in K^{\text{ur}}(p^n)$.\\
In this section we shall try to answer this question for an arbitrary formal group for $n = 1$ ie we would like to know if $\zeta_p \in K^{\text{ur}}(\pi)$.\\
First note that there is a formal group law $\mathbb{G}_m$ defined over $O_K$ given by $F(X,Y) = X + Y + XY$. This is a $\mathbb{Z}_p$ module and as $\mathbb{Z}_p$ module its height is 1. This implies $\mathbb{G}_m$ can not have structure of an unramified group law unless $K/\mathbb{Q}_p$ is unramified. It is easy to see, $\zeta_{p^n} - 1 \in \mathbb{G}_m[p^n]$ for all $n \geq 1$.\\~\\
We have the following lemma :\\~\\
\textbf{Lemma 4.1.1 :} Let $K$ be an unramified extension of $\mathbb{Q}_p$. Let $\mathfrak{F}$ be a formal group law defined over $O_K$. Then $\zeta_p \in K^{\text{ur}}(\mathfrak{F}[p])$. \\~\\
\textbf{Proof :} The $\mathbb{Z}_p$ modules $\mathbb{G}_m$ and $\mathfrak{F}$ both define unramified group laws over $O_K$. So the result follows from corollary-2.1.2 in \cite{several}. $\square$\\

Kolkata, India.\\
(First Version)\\
soumyadip.sahu00@gmail.com

\end{document}